\input amstex

\def\CC{{\Cal C}}
\def\CG{{\Cal G}}

\def\CM{{\Cal M}}

\def\CO{{\Cal O}}

\def\CT{{\Cal T}}
\def\CU{{\Cal U}}
\def\fg{{\frak g}}
\def\fm{{\frak m}}

\def\hO{{\hat{\Cal O}}}


\def\btu{\bigtriangleup}

\def\lra{\longrightarrow}

\parskip=6pt

\documentstyle{amsppt}
\document
\NoBlackBoxes
\nologo

\centerline{\bf LOCAL STRUCTURE OF MODULI SPACES}
\medskip
\centerline{\bf 3. Mathematische Arbeitstagung (Neue Serie)}
\smallskip
\centerline{\bf 20.-27 Juni 1997}
\bigskip
\centerline{Vadim Schechtman}
\smallskip
\centerline{Max-Planck-Institut f\"ur Mathematik,}
\centerline{Gottfried-Claren-Stra\ss e 26, 53225 Bonn}
\bigskip

\bigskip

This is a report on a joint work with Vladimir Hinich. 

{\bf 1.} Let $X$ be a smooth algebraic variety over a 
field $k$ of characteristic 
$0$. Assume that $X$ has no infinitesimal automorphisms, 
i.e. 
$H^0(X;\CT_X)=0$ ($\CT_X$ being the sheaf of 
vector fields). Let $\CM=Spec(\hO_{\CM})$ be the formal 
moduli space of deformations of $X$; $\hO_{\CM}$ is a 
complete local $k$-algebra. 
 
We have the {\it Kodaira-Spencer isomorphism}
$$
T_{\CM;X}=H^1(X;\CT_X)
\eqno{(KS)}
$$
Here $T_{\CM;X}=(\fm_{\CO_{\CM}}/\fm^2_{\CO_{\CM}})^*$ 
is the tangent space of $\CM$ at $X$. 
We want to consider the 

{\bf Problem.} {\it Describe the whole algebra $\hO_{\CM}$ 
in terms of $X$.} 

Let 
$$
\CT^{\bullet}_X=R\Gamma(X;\CT_X):\ 0\lra\CT^0\lra\CT^1\lra\ldots
$$
be a complex computing the sheaf cohomology of $\CT_X$. The last sheaf is a 
sheaf of Lie algebras, hence $\CT^{\bullet}$ may be chosen to be a {\it differential 
graded Lie algebra}; it is a correctly defined object of the appropriate 
derived category of {\it Homotopy Lie Algebras}, [HS1].

{\bf Theorem 1.} {\it One has a canonical isomorphism of $k$-algebras
$$
\hO_{\CM}=[H_0^{Lie}(\CT^{\bullet}_X)]^*
\eqno{(1)}
$$}

The homology of a (dg) Lie algebra is a (dg) coalgebra. The dual space is 
an algebra. The isomorphism (1) is a generalization of the Kodaira-Spencer 
isomorphism. 

This theorem is a just an example of a quite general fact; the similar 
results (with the same proof) hold true for other deformation problems. 
For example, we may 
wish to describe deformations of group representations, etc. Cf. [S2], [HS1].  

I know two proofs of Theorem 1. The first one works in the case when $\CM$ is 
smooth, and uses the {\it higher Kodaira-Spencer maps}, 
cf. [HS1]. 
The second one works in general situation. It uses certain 
very natural {\it sheaf property} of {\it Lie-Deligne functor}, and is decribed below. 

{\bf 2. Deligne groupoids.} Let $\fg=\oplus_{i\geq 0}\ \fg^i$ be a 
nilpotent dg Lie algebra. Recall that {\it groupoid} is a category 
with all morphisms being isomorphisms. The {\it Deligne groupoid} 
$\CG(\fg^{\bullet})$ is defined as follows. Its objects are {\it 
Maurer-Cartan elements} 
$$
MC(\fg^{\bullet}):=\{y\in \fg^1|dy+\frac{1}{2}[y,y]=0\}
$$
Let $\CG(\fg^0)$ be the Lie group corresponding to the 
nilpotent Lie algebra $\fg^0$. The algebra $\fg^0$ acts 
on $MC(\fg^{\bullet})$ by the rule
$$
x\circ y=dx+[x,y],\ \ \ x\in\fg^0, y\in MC(\fg^{\bullet}), 
$$
hence the group $\CG(\fg^0)$ acts on $MC(\fg^{\bullet})$. 
By definition, 
$$
Hom_{\CG(\fg^{\bullet})}(y,y')=\{g\in\CG(\fg^0)|y'=gy\}
$$
Morphisms are composed in the obvious way. Of course, 
this is a generalization of the Lie functor from Lie 
algebras to Lie groups. 

{\bf 3.} Let us return to our deformation situation. The 
variety $X$ defines a functor
$$
Def_X:\ Art_k\lra Groupoids
$$
where $Art_k$ is the category of artinian $k$-algebras 
with residue field $k$. Namely, $Def_X(A)$ is the groupoid 
whose objects are flat deformations of $X$ over $A$, 
and morphisms are isomorpshisms identical on $X$. 

On the other hand, if $\fg^{\bullet}$ is a dg Lie algebra 
over $k$, it defines a functor 
$$
\CG_{\fg^{\bullet}}:\ Art_X\lra Groupoids, 
$$
by
$$
\CG_{\fg^{\bullet}}(A)=\CG(\fm_A\otimes\fg^{\bullet}) 
$$
where $\fm_A$ is the maximal ideal of $A$. 

{\bf Example.} Assume that $X=Spec(R)$ is affine. Then 
one sees immediately from the definitions (Grothendieck) 
that one has an isomorphism of functors
$$
Def_X=\CG_{\fg^{\bullet}}
\eqno{(2)}
$$
where $\fg^{\bullet}=T_X=H^0(X;\CT_X)=Der_k(R)$ 
considered as a 
dg Lie algebra concentrated in dimension $0$. 

Sometimes when (2) holds, the people say that the dg 
Lie algebra $\fg^{\bullet}$ {\it governs the deformations 
of $X$}. 

{\bf Theorem 2.} {\it Let $X$ be arbitrary, and (2) holds 
for some $\fg^{\bullet}$. If $H^0(\fg^{\bullet})=0$ then 
$$
\hO_{\CM}=(H_0^{Lie}(\fg^{\bullet}))^*
$$}

{\bf Proof.} For an arbitrary $A\in Art_K$, we have
$$
Hom_{Art_k}(\hO_{\CM},A)=\pi_0(Def_X(A))=
\pi_0(\CG(\fm_A\otimes\fg^{\bullet}))=
$$
$$
=Hom_{alg}((H_0^{Lie}(\fg^{\bullet}))^*,A)\ \ \btu 
$$

{\bf 4.} Now, we know the Lie algebra $\fg^{\bullet}$ 
for affine varieties, and we want to know it for 
arbitrary ones, i.e. we want to glue them. 

Let $\fg^{\bullet}$ be a sheaf of nilpotent Lie algebras 
on a topological space $X$. It defines a presheaf 
of groupoids $\CG(\fg^{\bullet})$, 
$$
U\mapsto\ \CG(\Gamma(U;\fg^{\bullet})). 
$$

{\bf Theorem 3}, [H]. {\it Assume that 
$\fg^{\bullet}$ is a sheaf in the homotopy sense, i.e. for 
an open covering\ $\CU=\{U_i\},\ U=\bigcup U_i$, the natural map
$$
\Gamma(U;\fg^{\bullet})\lra
\check{\CC}(\CU;\fg^{\bullet})
$$
is quasiisomorphism. Then  $\CG(\fg^{\bullet})$ is a sheaf 
(i.e. stack).}

Here $\check{\CC}(\CU;\fg^{\bullet})$ is the \v Cech 
complex of the covering $\CU$.  

Theorem 1 follows immediately from Theorems 2 and 3, 
applied to a homotopy sheaf resolution of $\CT_X$.

{\bf 5.} It seems that the similar statements hold true 
in characteristic $p$, or in mixed characteristics, 
cf. [S2]. One should work with algebras of 
distributions instead of Lie algebras, and with 
simplicial objects instead of dg objects.

\bigskip
\centerline{\bf References}
\bigskip

[H] V.~Hinich, Descent of Deligne groupoids, {\it Int. Math. Res. Notices} 
{\bf 5}(1997), 223-239. 

[HS1] V.~Hinich, V.~Schechtman, Deformation theory and Lie algebra homology, 
alg-geom/9405013, {\it Alg. Coll.} (1997), to appear. 

[HS2] V.~Hinich, V.~Schechtman, On the descent of Deligne groupoids, 
Preprint (October 1994). 

[S1] V.~Schechtman, Localization of Deligne groupoids, Preprint (May 1994).

[S2] V.~Schechtman, Letter to B.Mazur (August 1994).

\enddocument